\documentstyle{amsppt}

\def\Prob{\text{\rm Prob}\,}
\topmatter
\title{Numerical characteristics of groups and corresponding relations}
\endtitle
\rightheadtext{Numerical characteristics of groups}
\author{A.~M.~Vershik}\endauthor
\affil{Steklov Institute of Mathematics}
{at St.Petersburg}
\endaffil
\date{Received June 27, 1999.}
\enddate
\thanks
{Partially supported by RFBR, grant 96-15-96060}
\endthanks
\endtopmatter

\document

\subhead\nofrills{1. Introduction}\endsubhead

Among various ways of ``measuring'' infinite groups, the most popular
is the so called growth, i.e. the asymptotics of the number of words of
given length. An enormous amount of works is devoted to this subject.
However, there are more deep characteristics, and their relations to the
growth and among themselves is maybe the most important subject of the
theory of numerical characteristics of algebraic systems. The most difficult
problem is to determine and calculate {\it invariant\/} characteristics
which do not depend on the choice of the system of generators. This field is
still little studied. In this short note we present related arguments
and examples which seem to be new. A more detailed presentation of these
results is to be exposed in author's paper in ``Uspekhi matematicheskih
nauk''.

We study mainly countable exponential groups (i.e. groups of
exponential
length) which are of the most interest for the applications we have in
mind. However, the general principles are valid for other cases as well. The
starting point for the author at one time was his studying the entropy of
random walks on groups. This notion (see below) was introduced in 1972 by
Avez \cite{1}, but I did not know this at the time and rediscovered it
several years later. However, the original purposes of introducing the
entropy were different from these of Avez. For the author, the entropy was
first of all an alternative characteristic with respect to the growth, and
the purpose was to study their relations. But very soon it became clear that the
entropy is a powerful technique for studying random walks, and first of all,
the Poisson boundaries which arose in many (not only probabilistic) problems
as well as in harmonic analysis on groups, etc. That is why the author
together with V.~A.~Kaimanovich (\cite{2, 3}) started studying random walks
from the entropic standpoint. Relations with amenability, Laplacian
spectrum, Kesten criterion, Furstenberg hypothesis etc. have united the entropy
with many other subjects. Subsequently this subject has intensely developed.
Here we return to the original purpose using already established relations and apply
the characteristics we introduce to new classes of groups. In particular, as
examples we continue studying so called local groups (\cite{6, 7}),
some of them (locally free groups) being recently considered in details in
a joint work of the author with S.~Nechaev and R.~Bikbov (\cite{9}).

In fact, the growth is related to only one of many measures on
the group, while the entropy enables to describe the asymptotics of
arbitrary measures. Relation between these characteristics is similar
to that between topological and metric entropy of transformations. The
principal fact -- {\it the fundamental inequality} -- is similar to the
inequality for these entropies. The analogy with dynamical systems is quite
far-reaching. Explicit calculation of the entropy is done in few
cases and it is by no means easy. The fundamental inequality itself
(see~below) in some cases was mentioned earlier in V.~Kaimanovich's work
(\cite{4}). In a recent paper (\cite{5}), a close question was considered
for Hausdorf dimensions of walks boundaries on trees. In terms of groups
this corresponds to the most studied case of free groups. However, it seems
that after this question was posed about twenty years ago, it has not been
studied in full generality.

\subhead\nofrills{2. Definitions}\endsubhead

Let $G$ be a countable group with a finite number of generators,
$({x_1,x_2,\ldots,x_k})$ be a system of generators, and
$(x_1^{-1},x_2^{-1},\ldots,x_k^{-1})$ be their inverses; denote by
$$
S=(x_1,x_2,\ldots, x_k,x_1^{-1}, x_2^{-1},\ldots,x_k^{-1}).
$$
Generally speaking, we do not assume that the system $S$ is minimal, i.e. it may be
that the system still generates the group after removing some elements.

Each element of the group is represented as an irreducible (i.e.
mutually inverse generators are not neighbours) words in the alphabet~$S$.
We distinguish words and elements of the group -- the length of a word is
considered literally, and the length of an element $g$, denoted by
$l_S(g)=l(g)$, i.e. the distance (in the words metric) from the group unity
in given generators, is the minimal length of the words representing this
element. The set of elements of length at most $n$ is denoted by
$W_{\leq n}$, and the set of elements of length exactly $n$ -- by $W_n$.
The number of
representations of a given element $g$ as a (maybe reducible) word of length
at most $n$ is denoted by $c_n(g)$. In terms of the Cayley graph, $W_{\leq n}$
is a ``ball'', i.e. the set of points which are at distance at most
$n$ from the unity, $W_n$ is a ``sphere'' of radius $n$, and
$c_n(g)$ is the number of paths of length at most $n$ from the unity $e$
to the element~$g$.

Let $m_{\leq n}$ ($m_n$ respectively) be the uniform normalized measure on
$W_{\leq n}$ (on $W_n$ respectively). Assume that $\mu_n$ is an arbitrary
symmetric ($\mu({g})=\mu({g^{-1}})$) measure on the set $S$.
Let $\mu^{*n}$ be the $n$th convolution of the measure $\mu$ with itself. It
is concentrated on the set $W_{\leq n}$. The convolution is defined by the
formula
$$
\mu^{*n}({g})=\sum_{g=g_1g_2,\ldots,g_n} \prod_{i=1,\ldots,n}\mu ({g_i)}.
$$
The uniform measure $\mu_S$ on $S$ is of particular interest. In this case
$$
\mu_S^{*n}({g})=c_n(g)/(2k)^n,
$$
thus the measure $\mu_S^{*n}({g})$ of an element $g$ is proportional to the
number of representations of this element in the system~$S$.

While the measures $m_{\leq n}$ and $m_n$ reflect in a sense the geometry of
a group, more exactly, the growth of vertices in the Cayley graph, the
measures $\mu_S^{*n}$ are related to dynamics, i.e. to the growth of
paths in this graph, and the measures
$\mu^{*n}$ -- to the growth of weighted number of paths. All these
measures depend on the choice of generators.

Now we can introduce a number of important constants. First of all, the well
known ``growth'' or volume.

\proclaim{Definition 1}
The logarithmic volume of a group $G$ in a given system of generators $S$
is the number
$$
\lim (\log |W_n|)/n =v.
$$
\endproclaim
The existence of the limit follows easily from the inequality
$$
|W_{n+m}| \leq |W_n| \cdot |W_m|.
$$

The logarithmic volume is not zero for groups of exponential growth and only
for them. Obviously, for these groups one can also calculate it by the
formula
$$
\lim \ \log (|W_{\leq n}|)/n =v,
$$
since the logarithmic asymptotics of balls and spheres are the same.

Recall that the entropy of a measure $\nu$ on a finite set $K$ is the number
$$
H(\nu)=-\sum_K \nu(k) \cdot \log (\nu(k))
$$
(logarithms are binary).

The entropy of a uniform measure is the logarithm of the number of
elements of the set.

\proclaim{Definition 2 {\rm(see~\cite {1, 3})}}
The entropy of a pair $(G, \mu)$ is the limit
$$
\lim  H(\mu^{*n})/n=h(G, \mu),
$$
i.e. the limit of normalized entropies of convolution powers of the measure
$\mu$ on the generators.
\endproclaim
This limit also exists and equals the infimum, as follows easily from the
inequality
$$
H(\nu_1*\nu_2) \leq H(\nu_1) + H(\nu_2).
$$
Note that the logarithmic volume is in fact the limit of normalized
entropies of uniform measures.

Let us introduce one more important characteristic. Consider the length of an
element $l(.)$ as a functional (depending on $S$) on the set
 $W_{\leq n}$ with the measure $\mu^{*n}$, and take its mean value with respect to
this measure. In other words, consider the expectation of the element length
with respect to the random walk for $n$ steps with initial measure
$\mu$, i.e. calculate the mean length of a random word $l_n(\mu)$.
There exists the limit
$$
l(\mu)= \lim l_n(\mu)/n,
$$
see Lemma~1.

\proclaim{Definition 3}
The number $l=l(\mu)$ is called the drift or the escape
with respect to a given measure, and in case of the measure
$\mu_S$ -- the escape
with respect to a given system of generators.
\endproclaim


\proclaim{Theorem 1. {\rm Fundamental inequality}}

The entropy of an arbitrary symmetric measure on a group concentrated on a
given system of generators of this group does not exceed the logarithmic
volume in this system multiplied by the escape
(drift) with respect to the measure,
$$
h \leq l \cdot v.
$$
\endproclaim

\subhead\nofrills{3. Proof of the theorem, corollaries, main problems
}\endsubhead

\proclaim{Lemma 1} If the drift is not zero, then the random sequence of
words lengths satisfies the Law of Large Numbers, i.e. for each
$\epsilon > 0$ there exists $N$ such that for all
$n \geq N$ the following relation holds,
 $$
 \mu^{*n}(g: |l(g)/l -1| \leq \epsilon) \geq (1-\epsilon).
 $$
\endproclaim

\demo{Proof}
Consider the lengths of initial segments of words as cylindric functionals
on the infinite product $(G^{\infty}, \mu^{\infty})$. Let $l_n(\omega)$
be the length of the element of the group corresponding to the initial
segment of $n$ letters of an infinite sequence $\omega$ in the alphabet $S$.
It is easy to see that the sequence
$l_n(\omega)$, $n=1,2,\dots$, satisfies the conditions of Kingman's Subadditive
Ergodic Theorem with respect to the left shift $T$ as a transformation of
$G^{\infty}$ with the product measure $\mu^{\infty}$, namely
$$
l_{n+m}(\omega) \leq l_n(\omega) + l_m(T^n \omega).
$$
Thus the sequence $l_n(\omega)/n$ converges almost everywhere. Its limit is
constant as follows from the $0$--$1$ law for a sequence of independent
variables and from the fact that this limit does not change under any change
of the initial segment of the sequence. The statement of the Lemma is
equivalent to the convergence in measure of the normalized lengths of
elements
and thus follows from the almost everywhere convergence.$\qquad\square$
\enddemo

\noindent
{\bf Remark.} If the drift is not zero, then it follows from the Lemma that the lengths of
almost all words grow linearly in $n$ with the same rate.
On the other hand, it is well known that for abelian and many other groups
the drift is zero and the rate of growth of the words lengths
is of order $\sqrt n$. The question posed in this connection is whether this
rate can for some groups be intermediate between $\sqrt n$ and $n$.
As far as we know, hitherto this question was not studied. The answer is
positive, the corresponding examples were constructed recently by a student
A.~Djubina on a basis of wreath products which are
quite useful in the theory of random walks and
were considered earlier
(see~\cite{8}). See in this volume
her result on properties of the drift on wreath products.

Now we pass to the proof of the Theorem.

\demo{Proof} Take an arbitrary positive $\epsilon$ and choose a positive
integer $n=n(\epsilon)$ satisfying the condition of Lemma~1. Decompose the measure $\mu^{*n}$ into a convex
combination of two measures, one of them $\mu_1$ being the normalized
restriction of this measure on the set
$V$ of elements whose lengths are in the interval
$[(1-\epsilon)l \cdot n, (1+\epsilon)l \cdot n]$, and the other
$\mu_2$ being the normalized restriction of the measure
$\mu^{*n}$on the complement of the set $V$ in
$W_{\leq n}$. By the Lemma, the measure of the set $V$ is greater than
$1-\epsilon$, and the number of elements of the second set does not exceed
the total number of elements in $W_{\leq n}$.

Then the entropy $H(\mu^{*n})$ may be estimated as
$$
H(\mu^{*n}) \leq (1-\epsilon)H(\mu_1)-\log (1-\epsilon)
                    +\epsilon \cdot \log (W_{\leq n})
$$
(the last summand estimates from above the entropy of the measure
$\mu_2$ by the entropy of the uniform measure on $W_{\leq n}$).
Divide by $n$ and take the limit on $n$, taking into account that the last
summand divided by $n$ tends to the logarithmic volume $v$ multiplied by
$\epsilon$ (see the definition). Thus we obtain
$$
h \leq (1-\epsilon) H(\mu_1)/n +\epsilon \cdot v.
$$

But $H(\mu_1)/n$ may be estimated from above by the expression
$(\log |V|)/n$ which gives $l \cdot v$ when $n$ tends to infinity, i.e.
$$
h \leq (1-\epsilon)l \cdot v + \epsilon \cdot v,
$$
and since $\epsilon$ is arbitrary, we obtain the desired inequality
$$
h \leq l \cdot v.
$$
$\square$
\enddemo

\proclaim{Corollary} If $l=0$, i.e. the growth of words is slower than
linear, then the entropy is zero, too. Thus the positiveness of the entropy
implies the linear growth of words typical with respect to this measure.
\endproclaim

Perhaps, the inverse statement is true, too\footnote{
V.~Kaimanovich has brought to my attention that the inverse statement
indeed follows from one estimation of N.~Varopoulos.}.

An analogy of the fundamental inequality with the inequality between
topological and metric entropy of dynamical systems (see, for
example,~\cite{11}) inevitably comes to mind.
It is convenient to write the inequality in the form
$$
h/l \leq v.
$$
The left hand side depends on the measure, and the right hand side depends
only on the set of generators (compare with~\cite{5}, where for trees the
right hand side is interpreted as the Hausdorf dimension of some measure). The
ratio in the left hand side is the
{\it normalized entropy} (i.e. the entropy on one step, or one letter, in
the real increment of the word). Denote it by~$\widehat h$.

If we consider the product measure on the infinite product of supports of
the measure $\mu$ with factors $\mu$, then the normalized entropy
$\widehat h$ may be easily interpreted as the limit entropy on one step of
partitions of infinite words on equivalence classes with respect to the
equality of initial segments of these words in the group, and $v$ becomes
the limit entropy of uniform measures on the classes. Thus our inequality
has not literally the same meaning as in topological dynamics (we have not a
transformation, but only a sequence of partitions).

However, this analogy is quite useful, first of all because of
similarity of main questions. Here they are.

Let the set of generators, and hence the value of $v$, be fixed. We call the
measure supported by a given set of generators and their inverses, that
provides the greatest $\widehat h$, the {\it measure of maximal normalized
entropy} with respect to a given set of generators.

\smallskip
A. WHEN FOR A GIVEN SET OF GENERATORS THERE EXISTS A SYMMETRIC MEASURE OF
MAXIMAL NORMALIZED ENTROPY? WHEN THIS MEASURE IS UNIQUE?

As the most natural measure $\mu$ one should take the uniform measure on
generators of the group and their inverses, but it seems that this measure does not
always give the maximal ratio $\widehat h=h/l$, since generally speaking the generators
may be not equal in rights.
\smallskip

It may be that the measure of maximal normalized entropy exists, but for a
given set of generators the maximal value of $\widehat h=h/l$ is strictly less
than~$v$. A symmetric set of generators in an exponential group is called
{\it extremal}, if there exists a symmetric measure on this set such that
for this measure the fundamental inequality in reduced form turns into
equality: $\widehat h =v$.

\smallskip
B. FOR WHICH GROUPS THERE EXISTS AN EXTREMAL SET OF GENERATORS?

\smallskip
Now we shall vary the set of generators too, and look for its invariant
characteristics.

Consider the expression
$$
\frac{h}{l \cdot v}=\frac{\widehat h}{v} \equiv q=q(S, \mu)
$$
(assuming that $l \ne 0$) and its supremum with respect to all symmetric
measures $\mu$ supported by a finite symmetric set $S$:
$$
q(S)= \sup_{\mu}  q(S,\mu).
$$
By Theorem~1, the number $q(S)$ does not exceed~$1$, and equals $1$ only
for extremal systems.
\smallskip
Thus we obtain an {\it objective method of comparing systems of generators:
the greater $q(S)$ is, the more effective the given system of generators
is}.

\smallskip
C. DESCRIBE INTRINSICALLY EXTREMAL SYSTEMS OF GENERATORS. CHOOSE THE ONE
WITH THE MAXIMAL VALUE OF $v$.
\smallskip

Note that the value of $v$ in itself does not tell much about the group, but
together with the following invariant it gives more deep information on the
group.

For an arbitrary exponential group, denote
$$
q_G= \sup_{S} q(S).
$$
This number is an invariant of the group. The system of generators which
provides the
maximal value of $q(S)$ we call {\it maximal}.
Groups with $q_G=1$ we call {\it free-like}. For the present one does not
know examples of groups with~$q_G<1$.

\smallskip
Let us elucidate the meaning of these notions. The greater the constant $q$
is, the greater the asymptotic fraction of elements of the group obtained by
the random walk (with a given measure on generators) for $n$ steps among all
elements of the corresponding length is. In particular, if $q=1$ for some
choice of generators and a proper measure on them (for example, uniform), then
typical random words exhaust asymptotically the whole group. In the opposite
case, they constitute only an exponentially small part.

\subhead\nofrills{4. Examples. Local groups}\endsubhead

Let us collect known facts on the fundamental inequality and examples of
calculating the introduced constants. There are quite little of them. Recall
that we consider only exponential groups. As to the growth (logarithmic
volume), there are dozens of works devoted to it. However, the entropy and
the drift are calculated in very few cases. The most effective way of
calculating the entropy is to use an analogue of Shannon's theorem proved
in~\cite{3} and~\cite{10}. We shall return to this calculation elsewhere.
Of course, it is interesting to do these calculations for classic groups
-- $SL(n,\bold Z)$, other matrix groups, solvable groups, etc. We
concentrate our attention on the class of so called local groups. But first
we give an example which is known for a long time.

\smallskip\noindent
{\bf Example 1. Free groups. }
Calculations in this case are well known.

\proclaim{Theorem 2}
The uniform measure on natural generators and their inverses in a free group with
$k$ ($k \geq 2$) generators has the maximal entropy, the fundamental
inequality turns into equality, and the system of generators is extremal,
i.e. a free group is free-like. The values of constants in this case are
the following:
$$
v= \log (2k-1),\; l=(k-1)/k,\; h=l \cdot v.
$$
\endproclaim

\demo{Proof}
Calculation of the logarithmic volume is trivial: the number of words of
length $n>0$ equals exactly
$2k \cdot (2k-1)^{n-1}$, hence $v=\log (2k-1)$. The drift, i.e. the
normalized on $n$ expectation of the word length, is the mean range of the
random walk on $\bold N$ with transition probabilities
$\Prob[i \to i+1]= 2k-1/2k$ and $\Prob[i \to i-1] = 1/2k-1$
for $i \ne 0$, hence $l=k-1/k$. The entropy of the measure $\mu^{*n}$
is equal to the entropy of the mixture of uniform distributions on words
of length $s=0,\ldots,n$ with binomial distribution as mixture
coefficients, and since the support of this distribution grows linearly and
not exponentially with $n$, its contribution to the entropy is zero.
Hence, the entropy is
equal to the entropy of the uniform distribution on words of the typical
length, i.e. $  \ln \cdot \log (2k-1) $, and after normalizing on $n$ we
obtain $h=k-1/k \cdot \log(2k-1)$, i.e. the equality in the fundamental
inequality.$\quad\square$
\enddemo

For extensions, factor groups obtained from free groups by subexponential
groups (so called virtually free groups) the situation is the same. However,
in case of general groups, the situation,  and in particular the calculation
of the entropy, which is the most
technically difficult part, is far more difficult.

The following question is initiated by similarity of hyperbolic and free
groups.

\smallskip
D.  ARE HYPERBOLIC IN SENSE OF GROMOV GROUPS FREE-LIKE?
\smallskip
If they are, then for which generators and measures the fundamental
inequality turns into equality? Also the following question is of interest:

\smallskip
E. ARE THE FUNDAMENTAL GROUPS OF SURFACES OF GENUS
$g \geq 2$ FREE-LIKE?
\medskip

\noindent
{\bf Example 2. Locally free groups, braid groups. }

We now recall the definition of local groups (\cite{6, 7}).
\proclaim{Definition 4}
A group (semi-group, algebra) is called local, if it is generated by a
finite or countable well-ordered set of generators $z_1,z_2,\ldots$,
the generators with numbers differing by more than $s$ (in the most
interesting case $s=1$) commute, and the neighbour ones may be connected by some
relations. If these relations depend only on the difference of numbers of
generators, then the local group (semi-group, algebra) is called
{\it local stationary}.
\endproclaim

We shall consider here only the case $s=1$, i.e. commutation at distance~$2$
and greater. The term ``local'' in algebra is evidently overloaded, but
here the analogy is rather physical -- commutation of remote elements
correspond to their ``non-interaction''.

Coxeter groups, braid groups, Hecke algebra etc. are examples of local
groups and algebras. More generally, one can consider as
$z_i$ not generators but subgroups or subalgebras. In the most recent time,
such objects become popular in mathematical physics (see~\cite{12, 13}).
Even more general notion is the localness with respect to a graph. In our
case the graph is a chain. Another interesting case is a cycle, when
generators are indexed with roots of unity.

\proclaim{Definition 5}
A locally free group $LF_k$ (semi-group $LF^{+}_k$, algebra)
is a group (semi-group, algebra respectively) such that the neighbour generators
$z_i, z_{i+1},  i=1 \dots k$,  $k \leq \infty$
generate free subgroups with two generators, and the other pairs of
generators commute.
\endproclaim

The locally free group $LF_k$ is in a sense an approximation of the braid
group $B_{k+1}$: it is simultaneously a subgroup and a supergroup of the braid
group.

The paper of S.~Nechaev, R.~Bikbov and the author~\cite{9} considers in
details the case of locally free group and semi-group from the above
standpoint. The results may be summed up as follows.

\proclaim{Theorem 3}
The natural generators in a locally free group (semi-group) are not
extremal. More exactly, for an arbitrary number of generators (including
infinite!) $k$ the following bounds hold.

For a locally free semi-group $LF^{+}_k$:
$$
  v = \log4 +  o(k^{-1}),
$$
$$
l=1,
$$
{\it since in the semi-group $LF^{+}_k$ there are no cancellations},
$$
   h=\log3+\epsilon(n,k),
$$
  where $\epsilon(n,k) \rightarrow 0$, $n,\,k \rightarrow \infty$,
   $k$ is the number of generators (maybe infinite), and $n$ is the length
of elements.

Thus for $k$ great enough and for infinite $k$ we obtain
$v \cdot l \not= h$.

For a locally free group $LF_k$:
$$
  v= \log 7 + o_k(1),
$$
$$
  l=2/3 - \beta,
$$
where $\beta$ is a non-negative constant (maybe zero for infinite $k$),
$$
h= \log (3-\alpha) + o(k^{-1}),
$$
where $\alpha$ is a constant which does not exceed~$1/2$ by absolute value.

Thus in natural generators of a locally free group (semi-group), the
fundamental inequality is strict, $v \cdot l > h$.
\endproclaim

The detailed calculations in~\cite{9} are related to combinatorics which is
interesting in itself (so called heaps and their bases) and is similar to the
combinatorics of Young diagrams and tableaus. We emphasize that for the
present only satisfactory bounds on constants are found
which enable us to deduce
non-maximality etc. But it seems that is is far more difficult
to find exact values of these constants.

\smallskip\noindent
{\bf Remarks.}

1. It is not known whether the uniform measure is the measure of maximal
entropy.

2. It is interesting that the logarithmic volume and the constants
 $l$ and $h$ in a locally free group {\it stabilize} with growth of the
number of generators $k$, and do not tend to infinity as one could suppose.

3. The paper~\cite{9} contains also calculations for locally free groups
with generators of finite order. If the order tends to infinity, all the
constants tend to their values for a locally free group.

4. From the Theorem and the above mentioned relation between the locally free
group and the braid group we obtain two-sided bounds on the same constants
for the braid group. It turns out that in this case, for the classical system
of generators of this group and the uniform measure, the fundamental inequality is not
strict either, and exactly as in a locally free group all constants stabilize to a finite
value when the number of generators tends to infinity.
Thus {\it in braid groups the set of elements which are typical random words
in natural generators with the uniform measure constitutes only an
exponentially small part of the set of all elements}.

Note also that generators in local groups are not equal in rights. More
exactly, a non-identity substitution which is not a reflection ($i \to k-i$)
does not generate an automorphism of the group. Thus one could
suppose that the uniform measure on generators and their inverses is not
quite natural. However, if we introduce a cycle locally free group (i.e.
identify the first and the last generators in a locally free group), then our
estimations will not change, though now a cycle group acts transitively on
generators and the uniform measure becomes natural.

\smallskip\noindent
 {\bf Example 3. Locally $M$-free groups. }
We now introduce the following class of local groups.
Consider some manifold of groups $M$ and a free group with two generators
$F_M(2)$ in this manifold. {\it Locally $M$-free group} $LF_M(k)$ with $k$
generators $(z_1, \dots, z_k)$ is the group generated by these generators
and the following relations:  $z_i \cdot z_j = z_j \cdot z_i, |i-j| \geq 2$,  and $(z_i, z_{i+1})$
are the generators of a subgroup which is canonically isomorphic to $F_M(2)$.
Of course, the existence of such groups for an arbitrary manifold $M$ is not
evident, and generally speaking the groups $LF_M(k)$ themselves do not belong to the
manifold for $k \geq 3$. For example, consider a locally free solvable group
of level~$2$ with two generators and construct a locally
solvable free group in this way. The structure of this group is rather complicated
(for $k \geq 3$) and it is already not solvable. However, it follows from the
above considerations of the locally free group that the growth in this group
stabilizes. These groups seem to arise in connection with crystal bases.

A more simple class of local groups is that of locally nilpotent free
groups. It arose in works on the theory of integrable systems. The simplest
variant of this group is the following. The neighbour generators $(z_i, z_{i+1})$
generate canonically the discrete Heisenberg group, and all the other pairs
commute. Though nilpotent groups have polynomial growth, locally nilpotent
free groups are already exponential for $k \geq 3$.
In all these cases the fundamental inequality remains presumably
strict. Relations with the Poisson boundary, whose non-triviality is
equivalent to positiveness of the entropy, with the spectral theory and harmonic
analysis on groups, in particular on local groups, are to be considered
elsewhere.

\smallskip
\eightpoint{
Vershik A. M. Numerical characteristics of the countable
groups and corresponding relations.

  The entropy of the random walk on the discrete countable group could
  be used for comparison of the system of the generators. Fundamental
  inequality between growth, entropy and escape gives the possibility
  to define ``the best'' system of the generators. We formulate a new
  circle of the problems related to the various growth and asymptotics
  on the groups.
}
\vskip-0.5cm


\Refs

\ref
\no 1
\by A.~Avez
\paper Entropie des groupes de type fini
\jour C. R. Acad. Sci. Paris
\vol 275A
\yr 1972
\pages 1363--1366
\endref

\ref
\no 2
\by A.~Vershik, V.~Kaimanovich
\paper Random walks on groups: boundaries, entropy, uniform distribution
\jour Sov. Dokl.
\vol 249 {\rm No.~1}
\yr 1979
\pages 15--18
\endref

\ref
\no 3
\by  A.~Kaimanovich, A.~Vershik
\paper Random Walks on Discrete Groups: Boundary and Entropy
\jour Ann. of Prob.
\vol 11, {\rm No.~3}
\yr 1983
\pages 457--490
\endref

\ref
\no 4
\by V.~Kaimanovich
\paper Entropy criterium of the maximality of boundary of walks on groups
\jour Sov. Dokl.
\vol 31
\pages 193--197
\yr 1985
\endref

\ref
\no 5
\by V.~Kaimanovich
\paper Hausdorf dimension of the harmonic measure on trees
\jour Ergod. Th. Dyn. Syst.
\vol 18 {\rm No.~3}
\pages 631--660
\yr 1998
\endref

\ref
\no 6
\by A.~Vershik
\paper Local stationary algebras
\jour Am. Math. Transl. (2)
\vol 148
\yr 1991
\pages 1--13
\endref

\ref
\no 7
\by A.~Vershik
\paper Local algebras and a new version of
Young's orthogonal form
\inbook Topics in Algebra. Banach Center Publ. {\bf 26}, part 2
\publ PWN-Polish Sci. Publ.
\publaddr Warszawa
\pages 467--473
\yr 1990
\endref

\ref
\no 8
\by A.~Vershik
\paper Amenability and approximation of infinite groups
\jour Selecta Math. Sov.
\vol 2, {\rm No.~4}
\yr 1982
\pages 311--330
\endref

\ref
\no 9
\by  R.~Bikbov, S.~Nechaev, A.~Vershik
\paper  Statistical properties
of braid groups in locally free approximation
\publ Preprint IHES (1999); submitted to Comm. Math. Phys.
\endref

\ref
\no 10
\by Y.~Deriennic
\paper Quelques applications du theoreme ergodic sous-additif
\jour Asterisque
\vol 74
\yr 1980
\pages 181--201
\endref

\ref
\no 11
\book Dynamical Systems, {\rm 2nd edition, Ya.~Sinai ed.}
\publ Springer
\yr 1999
\endref

\ref
\no 12
\by A.~Alekseev, L.~Faddeev, J.~Frohlich, and V.~Schmerus
\paper  Representation theory of lattice current algebras
\jour Commun. Math. Phys.
\vol 191
\pages 31
\yr 1998
\endref

\ref
\no 13
\by F.~Hausser, F.~Nill
\paper  Diagonal crossed product of
quantum groups
\jour Rev. Math. Phys.
\vol 11, {\rm No.~5}
\pages  553--630
\yr 1999
\endref

\endRefs
\enddocument